\documentclass[a4paper,twoside]{article} 

\usepackage{epsfig,amsmath,amsfonts}  

\newtheorem{lemma}{Lemma}

\newtheorem{prop}{Proposition}
\newtheorem{defin}{Definition}
\newtheorem{theo}{Theorem}

\newenvironment{proof}{{\em Proof:}}{\hfill\rule{2mm}{2mm}}

\newcommand{\C}{\mathbb C}  
  
\newcommand{\PP}{\mathbb P} 
\newcommand{\R}{\mathbb R}

\newcommand{\G}{{\cal G}}

\newcommand{\rank}{rank}
\newcommand{\grad}{grad}
\newcommand{\Perm}{Perm}
\newcommand{\im}{im}
\newcommand{\join}{\vee}
\newcommand{\meet}{\wedge}
\newcommand{\wtilde}{\widetilde}
\newcommand{\what}{\widehat}

\newcommand{\nodes}{\mbox{nodes}}

\newcommand{\ra}{\rightarrow}

\newcommand{\al}{\alpha}

\newcommand{\ga}{\gamma}

\newcommand{\ep}{\epsilon}
\newcommand{\la}{\lambda}
\newcommand{\om}{\omega}

\newcommand{\Ga}{\Gamma}
\newcommand{\De}{\Delta}

\newcommand{\La}{\Lambda}
\newcommand{\Om}{\Omega}

\newcommand{\Th}{\Theta}

\begin{document}

\title{Some applications of algebraic curves to computational vision
\thanks{This work is partially supported by the  Emmy Noether Institute for
Mathematics and the Minerva Foundation of Germany, by the Excellency Center
of the Israel Science Foundation "Group Theoretic Methods in the Study of
Algebraic Varieties" and by EAGER (European network
in Algebraic Geometry).}} 
\author{M. Fryers$^1$, J.Y. Kaminski$^2$, M. Teicher$^2$\\
$^1$ Department of Mathematics, \\
University of Hannover, Hannover, Germany. \\
$^2$ Department of Mathematics and Computer Science, \\
Bar-Ilan University, Ramat Gan, Israel.} 

\date{}
\maketitle

\section*{\centering Abstract}
{\em
We introduce a new formalism and a number of new results in the context of geometric computational vision. The classical scope of the research in geometric computer vision is essentially limited to static configurations of points and lines in $\PP^3$. By using some well known material from algebraic geometry, we open new branches to computational vision.    
 
We introduce algebraic curves embedded in $\PP^3$ as the building blocks from which the tensor of a couple of cameras (projections) can be computed. In the process we address dimensional issues and as a result establish the minimal number of algebraic curves required for the tensor variety to be discrete as a function of their degree and genus.   

We then establish new results on the reconstruction of an algebraic curves in $\PP^3$ from multiple projections on projective planes embedded in $\PP^3$. We address three different presentations of the curve: (i) definition by a set of equations, for which we show that for a generic configuration, two projections of a curve of degree $d$ defines a curve in $\PP^3$ with two irreducible components, one of degree $d$ and the other of degree $d(d-1)$, (ii) the dual presentation in the dual space $\PP^{3*}$, for which we derive a lower bound for the number of projections necessary for linear reconstruction as a function of the degree and the genus, and (iii) the presentation as an hypersurface of $\PP^5$, defined by the set of lines in $\PP^3$ meeting the curve, for which we also derive lower bounds for the number of projections necessary for linear reconstruction as a function of the degree (of the curve).    
 
Moreover we show that the latter representation yields a new and efficient algorithm for dealing with mixed configurations of static and moving points in $\PP^3$.  
}

\section{Introduction}
\label{sec:Intro}

Computational vision is the art of inferring three-dimensional properties of an object given several views of it. A view is concretely an image taken by a camera. From a mathematical point of view, a camera is a device that performs a projection from $\PP^3$ to $\PP^2$, through a center of projection. What are objects that we can deal with in this context? Until recently, for simplicity reasons, the objects that were considered in computer vision were essentially composed by static points and lines in $\PP^3$. In the last two decades, intensive research has been carried out in this context, analyzing the 3D structure of the scene when given multiple views. A summary of the past decade of work in this area can be found in \cite{Hartley-Zisserman-00,Faugeras-Luong-01} and references to earlier work in \cite{Faugeras-93}. 

The theory of multiple-view geometry and 3D reconstruction is somewhat fragmented when it comes to curve features, especially non-planar algebraic curves of any degree and genus. In \cite{Quan-96,Ma-Chen-94,Ma-Li-96} one shows how to recover the 3D position of a conic section from two and three views, and \cite{Schmid-Zisserman-98} shows how to recover the homography induced by a conic between two views of it, and \cite{Cross-Zisserman-98,Shashua-Toelg-ijcv} show how to recover a quadric surface from projections of its occluding conics. 

In all the above, the projection matrices are given but hardly any thing was done when this data is not available. \cite{Kahl-Heyden-98,Kaminski-Shashua-00} show how to recover the fundamental matrix from matching conics with the result that 4 matching conics are minimally sufficient for a finite numbers of solutions. \cite{Kaminski-Shashua-00} generalizes this result to higher order curves, but consider only planar curves.

In this paper we present a general theory for dealing with algebraic curves in computational vision, and more precisely in the multiple-view geometry and 3D reconstruction problem. Thus we investigate three of the fundamental questions of computer vision, using algebraic curves: (i) recovering camera geometry (fundamental matrix) from two views of the same algebraic curve, (ii) reconstruction of the curve from its projection across two or more views, and (iii) application of algebraic curves to the structure recovery of dynamic scenes from streams of projections (video sequences).

We start by giving a short presentation of the computer vision background related to our work. Then we show how one can recover the so-called projective stratum of the camera geometry (epipolar geometry) from two views of a (non-planar) algebraic curve. For that purpose we define and derive the {\it generalized Kruppa's equations} which are responsible for describing the epipolar constraint of two projections of the same algebraic curve. As part of the derivation of those constraints we address the issue of dimension analysis  and as a result establish the minimal number of algebraic curves required for the epipolar geometry to be defined up to a finite-fold ambiguity as a function of their degree and genus.

On the reconstruction front, the curve admits three different representations: (i) as the solution of a set of equations in $\PP^3$, for which we show that in a generic configuration, the reconstruction from two views of a curve of degree $d$ has two irreducible components, one of degree $d$ and the other of degree $d(d-1)$, (ii) as an hypersurface in the dual space $\PP^{3*}$, for which we derive a lower bound of the minimal number of projections necessary for linear reconstruction as a function of the curve degree and genus, and (iii) as an hypersurface in $\PP^5$, defined by the set of lines in $\PP^3$ meeting the curve, for which we also derive lower bounds for the number of projections necessary for linear reconstruction as a function of curve degree alone. Moreover we show that the latter representation yields a new and efficient algorithm for dealing with mixed configurations of static and moving points in $\PP^3$.  

\section{Foundation of linear computer vision}

Projective algebraic geometry provides a natural framework to geometric computer vision. However one has to keep in mind that the geometric entities to be considered are in fact embedded into the physical three-dimension euclidian space. For this content, the euclidian space is provided with three structures defined by three groups of transformations: the orthogonal group $\PP O_3$ (which defines the euclidian structure and which is included into the affine group), $\PP A_3$ (defining the affine structure and itself included into the projective group), $\PP G_3$ (defining the projective structure). We fix $[X,Y,Z,T]^T$, as homogeneous coordinates, and $T=0$ as the plane at infinity.

\subsection{A single camera system}
\label{SingleCamera}

Computational vision starts with images captured by cameras. The camera components are the following:
\begin{itemize}
\item a plane ${\cal R}$, called the {\it retinal plane}, or {\it image plane},
\item a point ${\bf O}$, called either the {\it optical center} or the {\it camera center}, and  which does not lie on the plane ${\cal R}$.
\end{itemize}

The plane ${\cal R}$ is regarded as a two dimension projective space embedded into $\PP^3$. Hence it is also denoted by $i(\PP^2)$. The camera is a projection machine:  $\pi: \PP^3 \setminus \{{\bf O}\} \rightarrow i(\PP^2), {\bf P} \mapsto \overline{{\bf O} {\bf P}} \cap i(\PP^2)$. The projection $\pi$ is determined (up to a scalar) by a $3 \times 4$ matrix ${\bf M}$ (the image of ${\bf P}$ being $\la {\bf M P}$). 

The physical properties of a camera imply that ${\bf M}$ can be decomposed as follows:
$$
{\bf M} = \left[\begin{array}{ccc}
	f & s & u_0 \\
	0 & \al f & v_0 \\
	0 & 0 & 1 
	\end{array}\right] [{\bf R}; {\bf t}],
$$
where $(f,\al,s,u_0,v_0)$ are the so-called internal parameters of the camera, whereas the rotation ${\bf R}$ and the translation ${\bf t}$ are the external parameters.

It is easy to see that:
\begin{itemize}
\item The camera center ${\bf O}$ is given by ${\bf MO}={\bf 0}$.
\item The matrix ${\bf M}^T$ maps a line in $i(\PP^2)$ to the only plane containing both the line and ${\bf O}$.
\item There exists a matrix $\what{\bf M} \in {\cal M}_{6 \times 3}(\R)$, which is a polynomial function of ${\bf M}$, that maps a point ${\bf p} \in i(\PP^2)$ to the line $\overline{{\bf Op}}$ (optical ray), represented by its Pl\"{u}cker coordinates in $\PP^5$. If the camera matrix is decomposed as follows:
$$
{\bf M}= \left[\begin{array}{c}
	\Ga^T \\
	\La^T \\
	\Th^T
\end{array}\right],
$$
then for ${\bf p}=[x,y,z]^T$, the optical ray ${\bf L}_{\bf p} = \what{{\bf M}}{\bf p}$ is given by the extensor: ${\bf L}_{{\bf p}}=x\La \meet \Th + y\Th \meet \Ga + z\Ga \meet \La$, where $\meet$ denotes the meet operator in the Grassman-Cayley algebra (see \cite{Barnabei-Brini-Rota-84}).
\item The matrix $\wtilde{\bf M}=\what{\bf M}^T$ maps lines in $\PP^3$ to lines in $i(\PP^2)$.
\end{itemize} 

Moreover we will need in the sequel to consider the projection of the {\bf absolute conic} onto the image plane. The absolute conic is simply defined by the following equations:
$$
\left\{\begin{array}{rcl}
X^2+Y^2+Z^2 & = & 0 \\
T & = & 0
\end{array}\right.
$$

By definition, the absolute conic is invariant by euclidian transformations. Therefore its projection onto the image plane, defined by the matrix $\om$, is a function of the internal parameters only. By Cholesky decomposition $\om={\bf L} {\bf U}$, where ${\bf L}$ (respectively ${\bf U}$) is lower (respectively upper) triangular matrix. Hence it is easy to see that ${\bf U}=\overline{\bf M}^{-1}$, where $\overline{\bf M}$ is the $3 \times 3$ matrix of the internal parameters of ${\bf M}$.

\subsection{A system of two cameras}

Given two cameras, $({\bf O}_j, i_j(\PP^2))_{j=1,2}$ are their components where $i_1(\PP^2)$ and $i_2(\PP^2)$ are two generic projective planes embedded into $\PP^3$, and  ${\bf O}_1$ and ${\bf O}_2$ are two generic points in $\PP^3$ not lying on the above planes. As in \ref{SingleCamera}, let $\pi_j: \PP^3 \setminus \{ {\bf O}_j \} \ra i_j(\PP^2), {\bf P} \mapsto \overline{{\bf O}_j {\bf P}} \cap i_j(\PP^2)$ be the respective projections. The camera matrices are ${\bf M}_i,i=1,2$.

\subsubsection{Homography between two images of the same plane}

Consider the case where the two cameras are looking at the same plane in space, denoted by $\De$. Let 
$$
{\bf M}_i = \left[\begin{array}{c}
	\Ga_i^T \\
	\La_i^T \\
	\Th_i^T
\end{array}\right]
$$
be the camera matrices, decomposed as above.  Let ${\bf P}$ be a point lying on $\De$. We shall denote the projections of ${\bf P}$ by ${\bf p}_i=[x_i,y_i,z_i]^T \cong {\bf M}_i {\bf P}$, where $\cong$ means equality modulo multiplication by a non-zero scalar. 

The optical ray generated by ${\bf p}_1$ is given by ${\bf L}_{{\bf p}_1}=x_1\La_1 \meet \Th_1 + y_1\Th_1 \meet \Ga_1 + z_1\Ga_1 \meet \La_1$. Hence ${\bf P}={\bf L}_{{\bf p}_1} \meet \De = x_1\La_1 \meet \Th_1 \meet \De + y_1\Th_1 \meet \Ga_1 \meet \De + z_1\Ga_1 \meet \La_1 \meet \De$. Hence ${\bf p}_2 \cong {\bf M}_2 {\bf P}$ is given by the following expression: ${\bf p}_2 \cong {\bf H}_\De {\bf p}_1$ where:
$$
{\bf H}_\De = \left[ \begin{array}{ccc} 
     \Ga_2^T (\La_1 \meet \Th_1 \meet \De)  & \Ga_2^T (\Th_1 \meet \Ga_1 \meet \De) & \Ga_2^T (\Ga_1 \meet \La_1 \meet \De) \\
     \La_2^T (\La_1 \meet \Th_1 \meet \De)  & \La_2^T (\Th_1 \meet \Ga_1 \meet \De) & \La_2^T (\Ga_1 \meet \La_1 \meet \De) \\
     \Th_2^T (\La_1 \meet \Th_1 \meet \De)  & \Th_2^T (\Th_1 \meet \Ga_1 \meet \De) & \Th_2^T (\Ga_1 \meet \La_1 \meet \De) \\
\end{array}\right].
$$
This yield the expression of the collineation ${\bf H}_\De$ between two images of the same plane.

\begin{defin}
The previous collineation is called the homography between the two images, through the plane $\De$.
\end{defin}

\subsubsection{Epipolar geometry}

\begin{defin}
Let $({\bf O}_j, i_j(\PP^2), {\bf M}_j)_{j=1,2}$ being defined as before. Given a pair $({\bf p}_1,{\bf p}_2) \in i_1(\PP^1) \times i_2(\PP^2)$, we say that it is a pair of corresponding or matching points if there exists ${\bf P} \in \PP^3$ such that ${\bf p}_j=\pi_j({\bf P})$ for $j=1,2$. 
\end{defin}

Consider a point ${\bf p} \in i_1(\PP^2)$. Then ${\bf p}$ can be the image of any point lying on the fiber $\pi_1^{-1}({\bf p})$. Then the matching point in the second image must lie on $\pi_2(\pi_1^{-1}({\bf p}))$, which is, for a generic point ${\bf p}$, a line on the second image. Since $\pi_1$ and $\pi_2$ are both linear, there exists a matrix ${\bf F} \in {\cal M}_{3 \times 3}(\R)$, such that: $\xi({\bf p})=\pi_2(\pi_1^{-1}({\bf p}))$ is given by ${\bf Fp}$ for all but one point in the first image.  

\begin{defin}
The matrix ${\bf F}$ is called the {\it fundamental matrix}, where as the line ${\bf l}_{\bf p}={\bf Fp}$ is called the {\it epipolar line} of ${\bf p}$.
\end{defin} 

Let ${\bf e}_1=\overline{{\bf O}_1 {\bf O}_2} \cap i_1(\PP^2)$ and ${\bf e}_2=\overline{{\bf O}_1 {\bf O}_2} \cap i_2(\PP^2)$. Those two points are respectively called the {\it first} and the {\it second epipole}. It is easy to see that ${\bf Fe}_1={\bf 0}$, since $\pi_1^{-1}({\bf e}_1)=\overline{{\bf O}_1 {\bf O}_2}$ and $\pi_2(\overline{{\bf O}_1 {\bf O}_2})={\bf e}_2$. Observe that by symmetry ${\bf F}^T$ is the fundamental matrix of the reverse couple of images. Hence ${\bf F}^T {\bf e}_2={\bf 0}$. Since the only point in the first image that is mapped to zero by ${\bf F}$ is the first epipole, ${\bf F}$ has rank $2$. 

Now we want to deduce an expression of ${\bf F}$ as a function of the camera matrices. By the previous analysis, it is clear that: ${\bf F}=\wtilde{\bf M}_2 \what{\bf M}_1$. Moreover we have the followings properties:

\begin{prop}
For any plane $\De$, not passing through the camera centers, the following equalities hold:
\begin{enumerate}
\item $$ {\bf F} \cong [{\bf e}_2]_\times {\bf H}_\De, $$ where $[{\bf e}_2]_\times$ is the matrix associated with the cross-product as follows: for any vector ${\bf p}$, ${\bf e}_2 \times {\bf p}=[{\bf e}_2]_\times {\bf p}$. Hence we have:
$$
[{\bf e}_2]_\times = \left[\begin{array}{ccc}
	0            & - {{\bf e}_2}_3  &  {{\bf e}_2}_2 \\
	{{\bf e}_2}_3   & 0             &  -{{\bf e}_2}_1 \\
	-{{\bf e}_2}_2  & {{\bf e}_2}_1    & 0
\end{array}\right].$$	
In particular, we have: ${\bf F}=[{\bf e}_2]_\times {\bf H}_\infty$, where ${\bf H}_\infty$ is the homography between the two images through the plane at infinity.
\item 
\begin{eqnarray} \label{FH} 
{\bf H}_\De^T {\bf F} + {\bf F}^T {\bf H}_\De={\bf 0}.
\end{eqnarray}
\end{enumerate}
\end{prop}
\begin{proof}
The first equality is clear according to its geometric meaning. Given a point ${\bf p}$ in the first image, ${\bf Fp}$ is its epipolar line in the second image. The optical ray ${\bf L}_{\bf p}$ passing trough ${\bf p}$ meets the plane $\De$ in a point ${\bf Q}$, which projection in the second image is ${\bf H}_\De {\bf p}$. Hence the epipolar line must be ${\bf e}_2 \join {\bf H}_\De {\bf p}$. This gives the required equality. The second equality is simply deduced from the first one by a short calculation.
\end{proof}

\begin{prop}
For a generic plane $\De$, the following equality holds:
$$
{\bf H}_\De {\bf e}_1 \cong {\bf e}_2.
$$
\end{prop}
\begin{proof}
The image of the ${\bf e}_1$ by the homography must be the projection on the second image of the point defined as being the intersection of the optical ray generated by ${\bf e}_1$ and the plane $\De$. Hence: ${\bf H}_\De {\bf e}_1 = {\bf M}_2 ({\bf L}_{{\bf e}_1} \meet \De)$. But ${\bf L}_{{\bf e}_1} = \overline{{\bf O}_1 {\bf O}_2}$. Thus the result must be ${\bf M}_2 {\bf O}_1$ (except when the plane is passing through ${\bf O}_2$), that is the second epipole ${\bf e}_2$.
\end{proof}

\subsubsection{Canonical stratification of the reconstruction}

Three-dimension reconstruction can be achieved from a system of two cameras, once the camera matrices are known. However a typical situation is that the camera matrices are unknown. Then we face a double problem: recovering the camera matrices and the actual object. There exists an inherent ambiguity. Consider a pair of camera matrices $({\bf M}_1, {\bf M}_2)$. If you change the world coordinate system by a transformation ${\bf V} \in \PP G_3$, the camera matrices are mapped to $({\bf M}_1 {\bf V}^{-1},{\bf M}_2 {\bf V}^{-1})$. Therefore we define the following equivalence relation:

\begin{defin}
Given a group of transformation $G$, two pairs of camera matrices, say $({\bf M}_1,{\bf M}_2)$ and $({\bf N}_1, {\bf N}_2)$ are said to be equivalent modulo multiplication by $G$ if there exists ${\bf V} \in G$, such that: ${\bf M}_1 =  {\bf N}_1 {\bf V} $ and ${\bf M}_2 = {\bf N}_2 {\bf V}$.
\end{defin}

Any reconstruction algorithm will always yields a reconstruction modulo multiplication by a certain group of transformations. More presicely there exist three levels of reconstruction according to the information that can be extracted from the two images and from a-priori knowledge of the world. \\

{\em Projective Stratum}

When the only available information is the fundamental matrix, then the reconstruction can be performed modulo multiplication by $\PP G_3$. Indeed from ${\bf F}$, the so-called intrinsic homography ${\bf S}=-\frac{{\bf e}_2}{\parallel {\bf e}_2 \parallel} {\bf F}$ is computed and the camera matrices are equivalent to: $([{\bf I};{\bf 0}],[{\bf S};{\bf e}_2])$. \\

{\em Affine Stratum}

When, in addition of the epipolar geometry, the homography between the two images through the plane at infinity, denoted by ${\bf H}_\infty$,  can be computed, the reconstruction can be done modulo multiplication by the group of affine transformations. Then the two camera matrices are equivalent to: $([{\bf I};{\bf 0}],[{\bf H}_\infty;{\bf e}_2])$.\\

{\em Euclidian Stratum}

The euclidian stratum is obtained by the data of the projection of the absolute conic $\Om$ onto the image planes, which allows recovering the internal parameters of the cameras. Once the internal parameters of the cameras are known the relative motion between the cameras expressed by a rotation ${\bf R}$ and a translation ${\bf t}$ can be extracted from the fundamental matrix. However only the direction of ${\bf t}$ but not the norm can be recovered. Then the cameras matrices are equivalent, modulo multiplication by the group of similarity transformations, to: $(\overline{\bf M}_1 [{\bf I};{\bf 0}], \overline{\bf M}_2 [{\bf R};{\bf t}])$, where $\overline{\bf M}_1$ and $\overline{\bf M}_2$ are the matrices of the internal parameters. 

Note that the projection of the absolute conic on the image can be computed using some a-priori knowledge of the world. Moreover there exist famous equations linking $\om_1$ and $\om_2$ when the epipolar geometry is given. These are the so-called {\em Kruppa's equations}, defined in the following proposition.

\begin{prop}
\label{KrupEq}
The projections of the absolute conic onto two images are related as follows. There exists a scalar $\la$ such that:
$$
[{\bf e}_1]_\times^T \om_1^* [{\bf e}_1]_\times = \la {\bf F}^T \om_2^* {\bf F},
$$
where $[{\bf e}_1]_\times$ is the matrix representing the cross-product by ${\bf e}_1$ and $\om_i^*$ is the adjoint matrix of $\om_i$.
\end{prop}

Let ${\ep}_i$ be the tangents to $\pi_i(\Om)$ through ${\bf e}_i$. Kruppa's equations simply state that ${\ep}_1$ and ${\ep}_2$ are projectively isomorphic. 

\section{Recovering the epipolar geometry from curve correspondences} 
\label{sec:EpipolarGeo}

Recovering epipolar geometry from curve correspondences requires the establishment of an algebraic relation between the two image curves, involving the fundamental matrix. Hence such an algebraic relation may be regarded as an extension of Kruppa's equations. In their original form (see proposition \ref{KrupEq}), these equations have been introduced to compute the camera-intrinsic parameters from the projection of the absolute conic onto the two image planes \cite{Maybank-Faugeras-92}. However it is obvious that they still hold if one replaces the absolute conic by any conic that lies on a plane that does not meet any of the camera centers. In this form they can be used to recover the epipolar geometry from conic correspondences \cite{Kahl-Heyden-98,Kaminski-Shashua-00}. Furthermore it is possible to extend them to any planar algebraic curve \cite{Kaminski-Shashua-00}. Moreover a generalization for arbitrary algebraic spatial curves is possible and is a step toward the recovery of epipolar geometry from matching curves.

Therefore we shall prove and generalize Kruppa's equations to arbitrary smooth irreducible curve that cannot be embedded in a plane and whose degree $d \geq 2$ (the case of line is excluded since one cannot deduce constraints on the epipolar geometry from a pair of matching lines and the case of planar curve has already been treated in \cite{Kaminski-Shashua-00}). 

At this stage we recall a number of facts, that will be useful in the sequel, concerning the projection of a space curve onto a plane. 

We shall mention that all our theoretical results are true when the ground field is the field of complex numbers. Finally we shall consider only the real solutions.

\subsection{Single view of a space curve}
\label{subsec:SingleView}

Let ${\bf M}$ be the camera matrix, ${\bf O}$ the camera center, ${\cal R}$ the retinal plane (as in \ref{SingleCamera}). Let $X$ be a smooth irreducible embedded in $\PP^3$ and which degree $d \geq 2$.  Let $Y$ be the projection of $X$ by the camera. It is well known that: 
\begin{enumerate}
\item The curve $Y$ will always contain singularities. Furthermore for a generic position of the camera center, the only singularities of $Y$ will be nodes.
\item The {\it class} of a planar curve is defined to be the degree of its dual curve. Let $m$ be the class of $Y$. Then $m$ is constant for a generic position of the camera center. 
\item If $d$ and $g$ are the degree and the genus of $X$, they are respectively, for a generic position of ${\bf O}$, the degree and the genus of $Y$ and the Pl\"ucker formula yields:
$$
\begin{array}{c}
m=d(d-1)-2(\sharp\nodes), \\ 
g=\frac{(d-1)(d-2)}{2} - (\sharp\nodes), 
\end{array}
$$
where $\sharp\nodes$ denotes the number of nodes of $Y$. Hence the genus, the degree and the class are related by: 
$$
m=2d+2g-2.
$$
\end{enumerate}

\subsection{Generalized Kruppa's equations} 
\label{subsec:ExtKrupEq}

We are ready now to investigate the recovery of the epipolar geometry from matching curves.  Let ${\bf M}_i$, $i=1,2$, be the camera matrices. Let ${\bf F}$ and ${\bf e}_1$ be the fundamental matrix and the first epipole, ${\bf F}{\bf e}_1=0$. We will need to consider the two following mappings: ${\bf p} \stackrel{\ga}{\mapsto} {\bf e}_1 \join {\bf p}$ and ${\bf p} \stackrel{\xi}{\mapsto} {\bf Fp}$. Both are defined on the first image plane; $\gamma$ associates a point to its epipolar line in the first image, while $\xi$ sends it to its epipolar line in the second image. 

Let $Y_1$ and $Y_2$ be the image curves (projections of $X$ onto the image planes). Assume that they are defined by the polynomials $f_1$ and $f_2$. Let $Y_1^\star$ and $Y_2^\star$ denote the dual image curves, whose polynomials are respectively $\phi_1$ and $\phi_2$.  Roughly speaking, the generalized Kruppa's equations state that the sets of epipolar lines tangent to the curve in each image are projectively equivalent.  

\begin{theo} \label{theoKrupeq} {\bf Generalized Kruppa's equations}\\ 
For a generic position of the camera centers with respect to the curve in space, there exists a non-zero scalar $\lambda$, such that for all points ${\bf p}$ in the first image, the following equality holds: 
\begin{equation}\label{ExtKrupEq} 
\phi_2(\xi({\bf p})) = \lambda \phi_1(\gamma({\bf p}))
\end{equation} 
\end{theo} 

{\it Remark:} Observe that if $X$ is a conic and ${\bf C}_1$ and ${\bf C}_2$ the matrices that respectively represent $Y_1$ and $Y_2$, the generalized Kruppa's equations reduce to the classical Kruppa's equations, that is: $[{\bf e}_1]_\times^T {\bf C}_1^\star [{\bf e}_1]_\times \cong {\bf F}^T {\bf C}_2^\star {\bf F}$, where ${\bf C}_1^\star$ and ${\bf C}_2^\star$ are the adjoint matrices of ${\bf C}_1$ and ${\bf C}_2$.  

\begin{proof}
Let $\epsilon_i$ be the set of epipolar lines tangent to the curve in image $i$. We start by proving the following lemma.

\begin{lemma}
The two sets $\epsilon_1$ and $\epsilon_2$ are projectively equivalent. Furthermore for each corresponding pair of epipolar lines, $({\bf l},{\bf l}') \in \epsilon_1 \times \epsilon_2$, the multiplicity of ${\bf l}$ and ${\bf l}'$ as points of the dual curves $Y_1^\star$ and $Y_2^\star$ are the same.
\end{lemma}
\begin{proof}
Consider the following three pencils:
\begin{enumerate}
\item $\sigma({\bf L}) \cong \PP^1$, the pencil of planes containing the baseline, generated by the camera centers,
\item $\sigma({\bf e}_1) \cong \PP^1$, the pencil of epipolar lines through the first epipole,
\item $\sigma({\bf e}_2) \cong \PP^1$, the pencil of epipolar lines through the second epipole.
\end{enumerate}
Thus we have $\epsilon_i \subset \sigma({\bf e}_i)$. Moreover if $E$ is the set of plane in $\sigma({\bf L})$ tangent to the curve in space, there exist a one-to-one mapping from $E$ to each $\epsilon_i$. This mapping also leaves the multiplicities unchanged. This completes the lemma.
\end{proof}

This lemma implies that both sides of the equation~\ref{ExtKrupEq} define the same algebraic set, that the union of eipolar lines through ${\bf e}_1$ tangent to $Y_1$. Since $\phi_1$ and $\phi_2$, in the generic case, have same degree (as stated in \ref{subsec:SingleView}), each side of equation~\ref{ExtKrupEq} can be factorized into linear factors, satisfying the following:
$$
\begin{array}{ccc}
\phi_1(\gamma (x,y,z)) &  =  & \prod_i (\alpha_{1i} x + \alpha_{2i} y + \alpha_{3i} z)^{a_i} \\ 
\phi_2(\xi (x,y,z)) & = & \prod_i \lambda_i (\alpha_{1i} x + \alpha_{2i} y + \alpha_{3i} z)^{b_i},
\end{array} 
$$
where $\sum_i a_i=\sum_j b_j=m$. By the previous lemma, we must also have $a_i=b_i$ for $i$. 
\end{proof}
 
By eliminating the scalar $\lambda$ from the generalized Kruppa's equations (\ref{ExtKrupEq}) we obtain a set of
bi-homogeneous equations in ${\bf F}$ and ${\bf e}_1$. Hence they define a variety in $\PP^2 \times \PP^8$. This gives rise to an important question. How many of those equations are algebraically independent, or in other words what is the dimension of the set of solutions? This is the issue of the next section.

\subsection{Dimension of the set of solutions}

Let $\{E_i({\bf F}, {\bf e}_1)\}_i$ be the set of bi-homogeneous equations on ${\bf F}$ and ${\bf e}_1$, extracted from the generalized Kruppa's equations (\ref{ExtKrupEq}). Our first concern is to determine whether all solutions of equation (\ref{ExtKrupEq}) are admissible, that is whether they satisfy the usual constraint ${\bf F}{\bf e}_1=0$. Indeed we prove the following statement:

\begin{prop}
As long as there are at least 2 distinct lines through ${\bf e}_1$ tangent to $Y_1$, equation (\ref{ExtKrupEq}) implies that $\rank{\bf F}=2$ and ${\bf F}{\bf e}_1={\bf 0}$.
\end{prop}
\begin{proof}
The variety defined by $\phi_1(\gamma({\bf p}))$ is then a union of at least 2 distinct lines through ${\bf e}_1$.  If equation (\ref{ExtKrupEq}) holds, $\phi_2(\xi({\bf p}))$ must define the same variety.

There are 2 cases to exclude: If $\rank{\bf F}=3$, then the curve defined by $\phi_2(\xi({\bf p}))$ is projectively equivalent to the curve defined by $\phi_2$, which is $Y_1^\star$. In particular, it is irreducible.

If $\rank{\bf F}<2$ or $\rank{\bf F}=2$ and ${\bf F}{\bf e}_1\neq{\bf 0}$, then there is some ${\bf a}$, not a multiple of ${\bf e}_1$, such that ${\bf Fa}={\bf 0}$.  Then the variety defined by $\phi_2(\xi({\bf p}))$ is a union of lines through ${\bf a}$. In neither case can this variety contain two distinct lines through
${\bf e}_1$, so we must have $\rank{\bf F}=2$ and ${\bf F}{\bf e}_1={\bf 0}$.
\end{proof}

As a result, in a generic situation every solution of $\{E_i({\bf F}, {\bf e}_1)\}_i$ is admissible.  Let $V$ be the subvariety of $\PP^2 \times \PP^8 \times \PP^2$ defined by the equations $\{E_i({\bf F},{\bf e}_1)\}_i$ together with ${\bf F}{\bf e}_1={\bf 0}$ and ${\bf e}_2{}^T {\bf F}={\bf 0}^T$, where ${\bf e}_2$ is the second epipole.  We next compute the lower bound on the dimension of $V$, after which we would be ready for the calculation itself.

\begin{prop} \label{prop:lowerbound}
If $V$ is non-empty, the dimension of $V$ is at least $7-m$.
\end{prop}
\begin{proof}
Choose any line ${\bf l}$ in $\PP^2$ and restrict ${\bf e}_1$ to the affine piece $\PP^2\setminus{\bf l}$. Let $(x,y)$ be homogeneous coordinates on ${\bf l}$.  If ${\bf F}{\bf e}_1={\bf 0}$, the two sides of equation (\ref{ExtKrupEq}) are both unchanged by replacing ${\bf p}$ by ${\bf p}+\alpha{\bf e}_1$.  So equation (\ref{ExtKrupEq}) will hold for all ${\bf p}$ if it holds for all ${\bf p}\in{\bf l}$. Therefore equation (\ref{ExtKrupEq}) is equivalent to the equality of 2 homogeneous polynomials of degree $m$ in $x$ and $y$, which in turn is equivalent to the equality of $(m+1)$ coefficients.  After eliminating $\lambda$, we have $m$ algebraic conditions on $({\bf e}_1,{\bf F},{\bf e}_2)$ in addition to ${\bf Fe}_1={\bf 0}$, ${\bf e}_2{}^T{\bf F}={\bf 0}^T$.

The space of all epipolar geometries, that is, solutions to ${\bf F}{\bf e}_1={\bf 0}$, ${\bf e}_2{}^T{\bf F}={\bf 0}^T$, is irreducible of dimension 7.  Therefore, $V$ is at least $(7-m)$-dimensional.
\end{proof}

For the calculation of the dimension of $V$ we introduce some additional notations. 
Given a triplet $({\bf e}_1,{\bf F},{\bf e}_2) \in \PP^2\times\PP^8\times\PP^2$, let $\{ {\bf
q}_{1\alpha}({\bf e}_1) \}$ (respectively $\{ {\bf q}_{2\alpha}({\bf e}_2) \}$) be the tangency points of the epipolar lines through ${\bf e}_1$ (respectively ${\bf e}_2$) to the first (respectively second) image curve. Let ${\bf Q}_\alpha({\bf e}_1,{\bf e}_2)$ be the 3d points projected onto $\{ {\bf q}_{1\alpha}({\bf e}_1) \}$ and $\{{\bf q}_{2\alpha}({\bf e}_2) \}$. Let ${\bf L}$ be the baseline joining the two camera centers. We next provide  a sufficient condition for $V$ to be discrete.

\begin{prop}
For a generic position of the camera centers, the variety $V$ will be discrete if, for any point $({\bf e}_1,{\bf F},{\bf e}_2) \in V$, the union of ${\bf L}$ and the points ${\bf Q}_\alpha({\bf e}_1,{\bf e}_2)$ is not contained in any quadric surface.
\end{prop}
\begin{proof}
For generic camera positions, there will be $m$ distinct points $\{{\bf q}_{1\alpha}({\bf e}_1)\}$ and $\{{\bf q}_{2\alpha}({\bf e}_2)\}$, and we can regard ${\bf q}_{1\alpha}$, ${\bf q}_{2\alpha}$ locally as smooth functions of ${\bf e}_1$, ${\bf e}_2$.

We let $W$ be the affine variety in $\C^3\times\C^9\times\C^3$ defined by the same equations as $V$.  Let $\Theta=({\bf e}_1, {\bf F}, {\bf e}_2)$ be a point of $W$ corresponding to a non-isolated point of $V$. Then there is a tangent vector $\vartheta=({\bf v}, \Phi, {\bf v}')$ to $W$ at $\Theta$ with $\Phi$ not a multiple of ${\bf F}$.

If $\chi$ is a function on $W$, $\nabla_{\Theta,\vartheta}(\chi)$ will denote the derivative of $\chi$ in the direction defined by $\vartheta$ at $\Theta$.  For
$$
\chi_\alpha({\bf e}_1,{\bf F},{\bf e}_2)= {\bf q}_{2\alpha}({\bf e}_2)^T{\bf F}{\bf q}_{1\alpha}({\bf e}_1),
$$
the generalized Kruppa's equations imply that $\chi_\alpha$ vanishes identically on $W$, so its derivative must also vanish.  This yields
\begin{equation}\label{EqNablaQFQ}
\nabla_{\Theta,\vartheta}(\chi_\alpha) =
(\nabla_{\Theta,\vartheta}({\bf q}_{2\alpha}))^T{\bf F}{\bf
q}_{1\alpha} \quad+{\bf q}_{2\alpha}^T\Phi{\bf q}_{1\alpha}+ {\bf
q}_{2\alpha}^T{\bf F}(\nabla_{\Theta,\vartheta}({\bf q}_{1\alpha}))
 = 0.
\end{equation}
We shall prove that $\nabla_{\Theta,\vartheta}({\bf q}_{1\alpha})$ is in the linear span of ${\bf q}_{1\alpha}$ and ${\bf e}_1$.  (This means that when the epipole moves slightly, ${\bf q}_{1\alpha}$ moves along the epipolar line.)  Consider $\kappa(t)=f({\bf q}_{1\alpha}({\bf e}_1 + t{\bf v}))$, where $f$ is the polynomial defining the image curve $Y_1$.  Since ${\bf q}_{1\alpha}({\bf e}_1 + t{\bf v}) \in Y_1$, $\kappa\equiv 0$, so the derivative $\kappa'(0)=0$.  On the other hand, $\kappa'(0)=\nabla_{\Theta,\vartheta}(f({\bf q}_{1\alpha})) =\grad_{{\bf q}_{1\alpha}}(f)^T\nabla_{\Theta,\vartheta}({\bf q}_{1\alpha})$.

Thus we have $\grad_{{\bf q}_{1\alpha}}(f)^T\nabla_{\Theta,\vartheta}({\bf q}_{1\alpha})=0$. But also $\grad_{{\bf q}_{1\alpha}}(f)^T{\bf q}_{1\alpha}=0$ and $\grad_{{\bf q}_{1\alpha}}(f)^T{\bf e}_1=0$.  Since $\grad_{{\bf q}_{1\alpha}}(f)\neq{\bf O}$ (${\bf q}_{1\alpha}$ is not a singular point of the curve), this shows that $\nabla_{\Theta,\vartheta}({\bf q}_{1\alpha})$, ${\bf q}_{1\alpha}$, and ${\bf e}_1$ are linearly dependent.  ${\bf q}_{1\alpha}$ and ${\bf e}_1$ are linearly independent, so $\nabla_{\Theta,\vartheta}({\bf q}_{1\alpha})$ must be in their linear span.

We have that  ${\bf q}_{2\alpha}^T{\bf F}{\bf e}_1= {\bf q}_{2\alpha}^T{\bf Fq}_{1\alpha}=0$, so ${\bf q}_{2\alpha}^T{\bf F} \nabla_{\Theta,\vartheta}({\bf q}_{1\alpha})=0$: the third term of equation (\ref{EqNablaQFQ}) vanishes.  In a similar way, the first term of equation (\ref{EqNablaQFQ}) vanishes, leaving
$$
{\bf q}_{2\alpha}^T\Phi{\bf q}_{1\alpha}=0.
$$
The derivative of $\chi({\bf e}_1,{\bf F},{\bf e}_2)={\bf F}{\bf e}_1$ must also vanish, which  yields:
$$
{\bf e}_2{}^T\Phi{\bf e}_1=0.
$$   
From the first equality, we deduce that for every ${\bf Q}_\alpha$, we have:
$$
{\bf Q}^T_\alpha {\bf M}_2^T \Phi {\bf M}_1 {\bf Q}_\alpha=0.
$$
From the second equality, we deduce that every point ${\bf P}$ lying on the baseline must satisfy:
$$
{\bf P}^T {\bf M}_2^T \Phi {\bf M}_1 {\bf P}=0.
$$
The fact that $\Phi$ is not a multiple of $F$ implies that ${\bf M}_2^T\Phi{\bf M}_1 \neq 0$, so together these two last equations mean that the union ${\bf L} \cup \{{\bf Q}_\alpha\}$ lies on a quadric
surface.  Thus if there is no such quadric surface, every point in $V$ must be isolated.
\end{proof}

Observe that this result is consistent with the previous proposition, since there always exist a quadric surface containing a given line and six given points. However in general there is no quadric containing a
given line and seven given points. Therefore we can conclude with the following theorem.

\begin{theo}
For a generic position of the camera centers, the generalized Kruppa's equations define the epipolar geometry up to a finite-fold ambiguity if and only if $m \geq 7$.
\end{theo}

Since different curves in generic position give rise to independent equations, this result means that the sum of the classes of the image curves must be at least $7$ for $V$ to be a finite set.

\section{3D Reconstruction} 
\label{3DReconstruction}

We turn our attention to the problem of reconstructing an algebraic curve from two or more views, given known camera matrices (epipolar geometries are known).  The basic idea is to intersect together the cones defined by the camera centers and the image curves. However this intersection can be computed in three different spaces, giving rise to different algorithms and applications. 

We shall mention that in \cite{Forsyth} a scheme is proposed to reconstruct an algebraic curve from a single view by blowing-up the projection. This approach results in a spatial curve defined up to an unknown projective transformation. In fact the only computation this reconstruction allows is the recovery of the properties of the curve that are invariant to projective transformation. Moreover this reconstruction is valid for irreducible curves only. However reconstructing from two projections not only gives the projective properties of the curve, but also the relative depth of it with respect to others objects in the scene and furthermore the relative position between irreducible components.

\subsection{Reconstruction $\PP^3$}

Let the camera projection matrices be ${\bf M}_1$ and ${\bf M}_2$. Hence the two cones defined by the image curves and the camera centers are given by: $\Delta_1({\bf P})=f_1({\bf M}_1 {\bf P})$ and $\Delta_2({\bf P})=f_2({\bf M}_2 {\bf P})$. The reconstruction is defined as the curve whose equations are $\Delta_1=0$ and $\Delta_2=0$. This curves has two irreducible components as the following theorem states.

\begin{theo}
For a generic position of the camera centers, that is when no epipolar plane is tangent twice to the curve $X$, the curve defined by $\{\Delta_1=0,\Delta_2=0\}$ has two irreducible components. One has degree $d$ and is the actual solution of the reconstruction. The other one has degree $d(d-1)$.
\end{theo}
\begin{proof}
For a line ${\bf l}\subset\PP^3$, we write $\sigma({\bf l})$ for the pencil of planes containing ${\bf l}$.  For a point ${\bf p}\in\PP^2$, we write $\sigma({\bf p})$ for the pencil of lines through ${\bf p}$. There is a natural isomorphism between $\sigma({\bf e}_i)$, the epipolar lines in image $i$, and $\sigma({\bf L})$, the planes containing both camera centers.  Consider the following covers of $\PP^1$:
\begin{enumerate}
\item $X \stackrel{\eta}{\longrightarrow} \sigma({\bf L}) \cong \PP^1$, taking a point $x \in X$ to the epipolar plane that it defines with the camera centers.
\item $Y_1 \stackrel{\eta_1}{\longrightarrow} \sigma({\bf e}_1) \cong \sigma({\bf L}) \cong \PP^1$, taking a point $y \in Y_1$ to its epipolar line in the first image.
\item $Y_2 \stackrel{\eta_2}{\longrightarrow} \sigma({\bf e}_2) \cong \sigma({\bf L}) \cong \PP^1$, taking a point $y \in Y_2$ to its epipolar line in the second image.
\end{enumerate}
If $\rho_i$ is the projection $X \to Y_i$, then $\eta=\eta_i\rho_i$. Let ${\cal B}$ the union set of branch points of $\eta_1$ and $\eta_2$.  It is clear that the branch points of $\eta$ are included in ${\cal B}$.  Let $S=\PP^1 \setminus {\cal B}$, pick $t \in S$, and write $X_S=\eta^{-1}(S)$, $X_t=\eta^{-1}(t)$.  Let $\mu_{X_S}$ be the monodromy: $\pi_1(S,t) \longrightarrow \Perm(X_t)$, where $\Perm(Z)$ is the group of permutation of a finite set $Z$.  It is well known that the path-connected components of $X$ are in one-to-one correspondence with the orbits of the action of $\im(\mu_{X_S})$ on $X_t$.  Since $X$ is assumed to be irreducible, it has only one component and $\im(\mu_{X_S})$ acts transitively on $X_t$.  Then if $\im(\mu_{X_S})$ is generated by transpositions, this will imply that $\im(\mu_{X_S})=\Perm(X_t)$.  In order to show that $\im(\mu_{X_S})$ is actually generated by transpositions, consider a loop in $\PP^1$ centered at $t$, say $l_t$.  If $l_t$ does not go round any branch point, then $l_t$ is homotopic to the constant path in $S$ and then $\mu_{X_S}([l_t])=1$.  Now in ${\cal B}$, there are three types of branch points:
\begin{enumerate}
\item branch points that come from nodes of $Y_1$: these are not branch points of $\eta$,
\item branch points that come from nodes of $Y_2$: these are not branch points of $\eta$,
\item branch points that come from epipolar lines tangent either to $Y_1$ or to $Y_2$: these are genuine branch points of $\eta$.
\end{enumerate}   

If the loop $l_t$ goes round a point of the first two types, then it is still true that $\mu_{X_S}([l_t])=1$.  Now suppose that $l_t$ goes round a genuine branch point of $\eta$, say $b$ (and goes round no other points in $\cal B$).  By genericity, $b$ is a simple two-fold branch point, hence $\mu_{X_S}([l_t])$ is a transposition.  This shows that $\im(\mu_{X_S})$ is actually generated by transpositions and so $\im(\mu_{X_S})=\Perm(X_t)$.

Now consider $\tilde{X}$, the curve defined by $\{\Delta_1=0,\Delta_2=0\}$.  By Bezout's Theorem $\tilde{X}$ has degree $d^2$.  Let $\tilde{x} \in \tilde{X}$.  It is projected onto a point $y_i$ in $Y_i$, such that $\eta_1(y_1)=\eta_2(y_2)$.  Hence $\tilde{X} \cong Y_1 \times_{\PP^1} Y_2$; restricting to the inverse image of the set $S$, we have $\tilde{X}_S \cong X_S \times_S X_S$. We can therefore identify $\tilde{X}_t$ with $X_t \times X_t$.  The monodromy $\mu_{\tilde{X}_S}$ can then be given by $\mu_{\tilde{X}_S}(x,y)=(\mu_{X_S}(x),\mu_{X_S}(y))$.  Since $\im(\mu_{X_S})=\Perm(X_t)$, the action of $\im(\mu_{\tilde{X}_S})$ on $X_t \times X_t$ has two orbits, namely $\{(x,x)\} \cong X_t$ and $\{(x,y) | x \neq y\}$.  Hence $\tilde{X}$ has two irreducible components.  One has degree $d$ and is $X$, the other has degree $d^2-d=d(d-1)$.
\end{proof} 
  
This result provides an algorithm to find the right solution for the reconstruction in a generic configuration, except in the case of conics, where the two components of the reconstruction are both admissible.

\subsection{Reconstruction in the Dual Space}

Let $X^\star$ be the dual variety of $X$. Since $X$ is supposed not to be a line, the dual variety $X^\star$ must be a hypersurface of the dual space \cite{Harris-92}. Hence let ${\Upsilon}$ be a minimal degree polynomial that represents $X^\star$. Our first concern is to determine the degree of $\Upsilon$. 

\begin{prop}
The degree of $\Upsilon$ is $m$, that is, the common degree of the dual image curves. 
\end{prop}
\begin{proof}
Since $X^\star$ is a hypersurface of $\PP^{3\star}$, its degree is the number of points where a generic line in $\PP^{3\star}$ meets $X^\star$. By duality it is the number of planes in a generic pencil that are tangent to $X$. Hence it is the degree of the dual image curve. Another way to express the same fact is the observation that the dual image curve is the intersection of $X^\star$ with a generic plane in $\PP^{3\star}$. Note that this provides a simple proof that the degree of the dual image curve is constant for a generic position of the camera center. 
\end{proof}

For the reconstruction of $X^*$ from multiple views, we will need to consider the mapping from a line ${\bf l}$ of the image plane to the plane that it defines with the camera center. Let $\mu: {\bf l} \mapsto {\bf M}^T {\bf l}$ denote this mapping. There exists a link involving $\Upsilon$, $\mu$ and $\phi$, the polynomial of the dual image curve: $\Upsilon(\mu({\bf l}))=0$ whenever $\phi({\bf l})=0$.  Since these two polynomials have the same degree (because $\mu$ is linear) and $\phi$ is irreducible, there exist a scalar $\lambda$ such that  
$$
\Upsilon(\mu({\bf l}))=\lambda \phi({\bf l}),
$$
for all lines ${\bf l} \in \PP^{2\star}$. Eliminating $\lambda$, we get ${m+2 \choose m} - 1$ linear equations on $\Upsilon$.  Since the number of coefficients in $\Upsilon$ is ${m+3 \choose m}$, we can state the following result: 

\begin{prop} \label{dualNbViews}
The reconstruction in the dual space can be done linearly using at least $k \geq \frac{m^2+6m+11}{3(m+3)}$ views. 
\end{prop}

From a practical point of view, it is worth noting that the fitting of the dual image curve is not necessary. It is sufficient to extract tangents to the image curves at distinct points. Each tangent ${\bf l}$ contributes to one linear equation on $\Upsilon$: $\Upsilon(\mu({\bf l}))=0$. However one cannot obtain more than ${m+2 \choose m}-1$ linearly independent equations per view. \\ 

\subsection{Reconstruction in $\G(1,3)$}  

Let $\G(1,3)$ be the Grassmanian of lines of $\PP^3$. Consider the set of lines in $\PP^3$ meeting the curve $X$ of degree $d$. This defines an irreducible subvariety of $\G(1,3)$ which is the intersection of $\G(1,3)$ with an irreducible hypersurface of degree $d$ in $\PP^5$ (see \cite{Harris-92}), given by a homogeneous polynomial $\Gamma$, defined modulo the $d$th graded piece $I(\G(1,3))_d$ of the ideal of $\G(1,3)$ and modulo scalars. However picking one representative of this equivalence class is sufficient to reconstruct entirely without any ambiguity the curve $X$. Hence we need to compute the class of $\Ga$ in the homogeneous coordinate ring of $\G(1,3)$, or more precisely in its $d$th graded piece, $S(\G(1,3))_d$, which dimension is ${d+5 \choose d} - {d-2+5 \choose d-2}$. 

Let $f$ be the polynomial defining the image curve, $Y$. Consider the mapping  that associates to an image point its optical ray: $\nu: {\bf p} \mapsto \widehat{\bf M} {\bf p}$. Hence the polynomial $\Gamma(\nu({\bf p}))$ vanishes whenever $f({\bf p})$ does. Since they have same degree and $f$ is irreducible, there exists a scalar $\lambda$ such as for every point ${\bf p} \in \PP^2$, we have:   
$$
\Gamma(\nu({\bf p})) = \lambda f({\bf p}).
$$
This yields ${d+2 \choose d}-1$ linear equations on $\Gamma$. 

Hence a similar statement to that in Proposition \ref{dualNbViews} can be made:
\begin{prop}
The reconstruction in $\G(1,3)$ can be done linearly using at least $k \geq \frac{1}{6}\frac{d^3+8d^2+23d+28}{d+3}$ views. 
\end{prop}

As in the case of reconstruction in the dual space, it is not necessary to explicitly compute $f$. It is enough to pick points on the image curve. Each point yields a linear equation on $\Gamma$: $\Gamma(\nu({\bf p}))=0$. However for each view, one cannot extract more than $\frac{1}{2}d^2+\frac{3}{2}d$ independant linear equations.  

\section{Applications to dynamic configurations of points}

In this section we show that the reconstructions in $\G(1,3)$ can be applied to trajectory recovery. Consider a point moving along a smooth trajectory. The motion is assumed to be well approximated by a low degree irreducible algebraic curve. This requirement is in fact very natural and has a widely broad validity in practice. Now we proceed to show how trajectory recovery can be achieved.  

A set of cameras ${\bf M}_i, i=1,...,m$ which are either static or moving is viewing at a set of points ${\bf P}_j, j=1,...,n$ either static or moving. The cameras are independant and in particular they are not supposed to be synchronized. Let ${\bf p}_{ijk}$ be the projection of the point ${\bf P}_j$ onto the camera $i$ at time $k$.  

For a given point ${\bf P}_j$, for all $i$ and $k$, the optical rays, ${\bf L}_{ijk}=\what{\bf M}_i {\bf p}_{ijk}$ meet the trajectory of ${\bf P}_j$. Then according to the geometric entity those rays generate, the motion of ${\bf P}_j$ can be recovered. Here we provide a table that gives the correspondence between this entity and the motion of ${\bf P}_j$. 

\begin{center}
\begin{tabular}{|c|c|}
\hline
Motion of ${\bf P}_j$ & Geometry entity generated \\ 
                      & by $\{{\bf L}_{ijk}\}$ in $\PP^5$ \\ 
\hline
Static point & Plane included \\ 
             & in $\G(1,3)$ \\ 
\hline
Point moving on a line & Intersection of $\G(1,3)$ \\ 
                       & with an hyperplane\\ 
\hline
Point moving on a conic & Intersection of $\G(1,3)$ \\
                        & with a quadric\\ 
\hline
... & ... \\
\hline
Point moving on a curve & Intersection of $\G(1,3)$ \\
of degree $d$ &  with an hypersurface of degree $d$ \\ 
\hline
\end{tabular}
\end{center}

\end{document}